\documentstyle[11pt,amssymb,amstex]{amsart}

\newtheorem{thm}{Theorem}[section]
\newtheorem{lem}[thm]{Lemma}
\newtheorem{cor}[thm]{Corollary}
\newtheorem{prop}[thm]{Proposition}

\newtheorem{rem}[thm]{Remark}

\def\a{\alpha}

\begin{document}

\title[dissipative quadratic stochastic operators ]
{The Limit Behavior of Trajectories of Dissipative Quadratic Stochastic Operators on finite-dimensional simplex}

\maketitle
\begin{center}
\author{F.A.Shahidi $^{\rm a}$\footnote{Corresponding author. Email: farruh.shahidi@@gmail.com} and M.T. Abu Osman $^{\rm a}$}\\\vspace{12pt}
\address{$^{\rm a}${\em{Faculty of Information and Communication Technology, International Islamic University Malaysia, P.O Box 10, 50728, Kuala Lumpur, Malaysia}}}

\end{center}

\begin{abstract}

The limit behavior of trajectories of dissipative quadratic stochastic operators on a finite-dimensional simplex is fully studied. It is shown that any dissipative quadratic stochastic operator has either unique or infinitely many fixed points. If dissipative quadratic stochastic operator has a unique point, it is proven that the operator is regular at this fixed point. If it has infinitely many fixed points, then it is shown that $\omega-$ limit set of the trajectory is contained in the set of fixed points.\vskip 0.3cm \noindent {\it
Mathematics Subject Classification}: 15A51, 47H60, 46T05, 92B99.
\end{abstract}

\section{Introduction}
Many situations in population genetics can be described by discrete dynamical systems \cite{Ber},\cite{HS},\cite{K},\cite{Lyu}.
More precisely, dynamical systems in population genetics are generated by nonlinear maps from finite-dimensional space into itself. S.N.Bernstein in \cite{Ber} introduced a class of dynamical systems which is now called a \textit{quadratic stochastic operator} (q.s.o. in short). Q.s.o. has the following form:
$$
(Vx)_{k}=\sum\limits_{i,j=1}^{m}p_{ij,k}x_{i}x_{j},  $$
where $x=(x_1,x_2,...,x_m)\in S^{m-1}$ and

$$S^{m-1}=\{x\in R^{m}: x_{i}\geq 0, \ \
\sum\limits_{i=1}^{m}x_{i}=1\}$$
is $m-1$ dimensional standard simplex.
Here the coefficients $p_{ij,k},$ called  heredity coefficients, satisfy the following conditions
$$p_{ij,k}=p_{ji,k}\geq 0 , \ \ \sum\limits_{k=1}^mp_{ij,k}=1.$$

If the above conditions are satisfied, then one can easily check that q.s.o. maps the simplex into itself.

Research on quadratic stochastic operators was later developed by a number of authors (see for example \cite{Ga1}, \cite{GEsh}, \cite{RozZ}, \cite{Fsib}). The main problem in this area is the study of limit behavior of the trajectory of an initial point taken from the simplex. Note that the limit behavior of trajectories of quadratic stochastic operator was fully classified on one dimensional simplex by Yu.I.Lyubich in \cite{Lyu0}, \cite{Lyu}. However, the problem is still open in higher dimensions (even in two-dimensional simplex). In some papers, classes of q.s.o. were considered with additional properties that enables the study easier(see \cite{Roz}). For example, in \cite{Ga1}, the class of Volterra q.s.o. was outlined which satisfies the additional condition $p_{ij,k}=0$ whenever $k\notin \{i,j\}.$ It was shown that the trajectory of any nonfixed initial point, approaches the bound of the simplex.

A q.s.o. is called \textit{doubly stochastic} if $Vx\prec x$ for all $x\in S^{m-1},$ where $"\prec"$ is the notation of majorization, that is the classical ordering by comparison of partial sums of the coordinates after the nonincreasing rearrangement. The class of doubly stochastic operators was introduced in \cite{Ga2},\cite{GSh}, and the structural properties of the set of doubly stochastic operators were studied. The class of doubly stochastic quadratic operators is sufficiently large, therefore the study of limit behavior is complicated. However, by altering $Vx$ and $x$ in the definition of doubly stochastic operator, we find the  class of dissipative q.s.o. The class of dissipative q.s.o. was introduced in \cite{Sh1}, in which it was studied some properties of dissipative q.s.o.. Here we note that the definition of dissipativity is not given in terms of wandering set as in \cite{Nic}, but in terms of majorization.

From the definition of dissipative and doubly stochastic q.s.o. one can see that the intersection of the set of dissipative and doubly stochastic q.s.o.'s is nothing but permutation operators of the simplex. Besides, in \cite{Sh1}, it was shown that only identity operator can simultaneously be Volterra and dissipative. So, we assert that the class of dissipative operators is the new class of q.s.o., which makes the study of this class more interesting.

We say that a q.s.o. $V$ satisfies the \textit{ergodic theorem }(\cite{Ul}), if for any $x\in S^{m-1}$ the limit

$$
\lim\limits_{n\rightarrow\infty}\frac{x+Vx+\cdots+V^{n-1}x}{n}
$$
exists. In \cite{Sh1} it was shown that dissipative operator satisfies the ergodic theorem. However, the convergence of Cesaro mean of the trajectory of dissipative q.s.o. does not imply the convergence of the trajectory of the dissipative q.s.o. Therefore, the study of limit behavior of the trajectory of dissipative q.s.o. is still open.

In the present paper, we classify the limit behavior of dissipative q.s.o. on finite dimensional simplex. The paper is organized as follows. In section 2 we give some terminology and notations. Also, we provide some previous results concerning dissipative q.s.o. In section 3 we give main results of the paper.

\section{Preliminaries}

In this section we give some definitions and state some previous results.
Let $$S^{m-1}=\{x\in R^{m}: x_{i}\geq 0, \ \
\sum\limits_{i=1}^{m}x_{i}=1\}$$

be a $(m-1)-$dimensional simplex.
Then the vectors
$$e_k=(0,0,\cdots,\underbrace{1}_{k},\cdots,0),$$
($ k=\overline{1,m}$) are its vertices. For $\a\subset I={1,2,\cdots,m},$ the set
$F_{\a}=\{x\in S^{m-1}: x_i=0,\ i\notin\a \}$ is called a \textit{face} of the simplex.

For $x\in S^{m-1}$ Let us put
$x_{\downarrow}=(x_{[1]},x_{[2]},\cdots x_{[m]}),$ where
$(x_{[1]},x_{[2]},\cdots x_{[m]})$- nonincreasing rearrangement of
$(x_1,x_2,\cdots x_m),$ that is $x_{[1]}\geq x_{[2]}\geq \cdots
\geq x_{[m]}.$

\textbf{Definition 1}\cite{Ma}. We say that $x$ is {\it majorized} by $y$ on $S^{m-1}$, and write $x\prec y$(or $y\succ x$) if
:
$$ \sum\limits_{i=1}^{k}x_{[i]}\leq
\sum\limits_{i=1}^{k}y_{[i]}, \ \ k=\overline{1,m-1}.$$

It is easy to see that for any $x\in S^{m-1}$ we have
$$(\frac1m,\frac1m,\cdots,\frac1m)\prec x\prec (1,0,0,\cdots,0).$$

The following proposition can be found in \cite{Ma}.
\begin{prop}The following statements are mutually equivalent

1) $x\prec y,$ that is  $ \sum\limits_{i=1}^{k}x_{[i]}\leq
\sum\limits_{i=1}^{k}y_{[i]}, \ \ k=\overline{1,m-1}.$

2) $x=Py$ for some doubly stochastic matrix $P.$

3) Vector $x$ belongs to the convex hull of the $m!$
vectors obtained by permuting the components of the vector $y.$
\end{prop}

\textbf{Definition 2.} An operator  $V:S^{m-1}\rightarrow S^{m-1}$ is called {\it
dissipative} if
$$
Vx\succ x, \ \ \forall x\in S^{m-1}.\eqno (1)
$$

The reason that the operator with above property is called dissipative is simple. Since $Vx\succ x,$
then according to the above proposition the point $Vx$ does not belong to the convex hull of the points, obtained by permuting the components of $x.$ The same can be said about points $V^2x$ and $Vx$, $V^3x$ and $V^2x$, and so on. So, if we consider points $x, Vx, V^2x,\cdots$, then these points dissipate inside the simplex.

It was stated in \cite{Sh1} dissipative linear operators are just permutations.
Well, it is not generally so. Here the

\textbf{Observation.} Let us consider the case when $V$ is a linear
dissipative operator, that is $Vx=Ax$, where
$A=(a_{ij})_{i,j=\overline{1,m}}$ is an $m\times m$ matrix. Since
$Vx\succ x,$ then by putting $x=e_i$ we have $Ae_i \succ e_i.$
At the same time, it is easy to see that $Ae_i \prec e_i.$ That is why
$(Ae_i)_{\downarrow}=(e_i)_{\downarrow},$ which means that only
one component of the vector $Ae_i$ is $1$ and the others are $0$. Hence, the operator can
either be permutation operator or maps the simplex into its face which makes the
study of the operator easy in either way.
Here the example of linear dissipative operator which is not permutation:
$$V:(x_1,x_2,x_3)\rightarrow(x_1+x_2,x_3,0).$$
However, if we consider the interior of the simplex, then linear dissipative
operators are only permutations. Therefore, it is more interesting to study non-linear dissipative ones. In what follows, we consider quadratic operators.
Recall \cite{Lyu} a stochastic operator $V:S^{m-1}\rightarrow S^{m-1}$ is called a {\it quadratic stochastic
operator } if it has a form
$$
(Vx)_{k}=\sum\limits_{i,j=1}^{m}p_{ij,k}x_{i}x_{j}, \ \ \ \eqno
(2)
$$
where $x=(x_1,x_2,...,x_m)\in S^{m-1}$. Here the coefficients
$p_{ij,k}$ satisfy the following conditions
$$p_{ij,k}=p_{ji,k}\geq 0 , \ \ \sum\limits_{k=1}^mp_{ij,k}=1.\eqno (3)$$
More information concerning q.s.o. can be found in \cite{Lyu}. We note that linear operator on
$S^{m-1}$  can be considered as a quadratic operator since we can always multiply the operator by $\sum\limits_{i=1}^mx_i=1$. But, since the study of the limit behavior trajectories  of linear dissipative operators are simple, here and henceforth we exclude this kind of operators from the set of dissipative q.s.o.

Let us consider the following operator
\begin{eqnarray*}
&&(Vx)_1=x_1^2+x_2^2+x_3^2+x_1x_2+x_1x_3+x_2x_3,\\
&&(Vx)_2=x_1x_2+x_1x_3,\\
&&(Vx)_3=x_2x_3.
\end{eqnarray*}

Simple calculations show that $(Vx)_1\ge max\{x_1,x_2,x_3\}$ and
$(Vx)_1+(Vx)_2\ge max\{x_1+x_2,x_1+x_3\}, \ (Vx)_1+(Vx)_3\ge x_2+x_3.$
Therefore, above operator is dissipative quadratic stochastic operator.

From the definitions 1 and 2, it is easy to see that if an operator is dissipative,
then by changing its components it preserves its dissipativity. This gives more examples of dissipative q.s.o.

Now, for the purpose of completeness we state two Lemmas from \cite{Sh1} and provide their short proofs
since we are going to use them frequently.

Given q.s.o. $V$ we denote $a_{ij}=(p_{ij,1},p_{ij,2},\cdots
p_{ij,m})\ \ \forall i,j=\overline {1,m}$, where $p_{ij,k}$ are
the coefficients of q.s.o. $V$. One can see that
$a_{ij}\in S^{m-1}$, for all $i,j=\overline{1,m}.$

\begin{lem}\label{l1} Let $V$ be a dissipative q.s.o.
Then the following conditions hold
$$ (a_{ii})_{\downarrow} =e_1 \ \ \forall i=\overline{1,m}.$$
\end{lem}

\begin{pf} Due to dissipativity of $V$ one has $Vx\succ x, \ \ \forall x\in S^{m-1}.$
Now by putting $x=e_i$ we get $e_i\prec Ve_i.$ On the other hand,
we have $e_i\succ x, \ \ \forall x\in
S^{m-1}.$ That is why $(Ve_i)_{\downarrow}=(e_i)_{\downarrow}=e_1.$
Then the equality $Ve_i=a_{ii}$ implies the assertion.
\end{pf}

The above lemma implies that any dissipative q.s.o. can be written as
$$(Vx)_k=\sum\limits_{i\in
\a_k}x_i^2+2\sum\limits_{i<j}p_{ij,k}x_ix_j \ \
k=\overline{1,m},\eqno (4)$$ where

$$\a_k\subset I=\{1,2\cdots,m\},\ \a_i\cap \a_j=\emptyset,\ i\neq j,\
\bigcup\limits_{k=1}^m\a_k=I. \eqno (4')$$

\begin{lem}\label{l2} Let (4) be a dissipative q.s.o.
\begin{itemize}
\item[(i)] If $j\in \alpha_{k_0},$ then
$p_{ij,k_0}=(a_{ij})_{[1]}\geq \frac{1}{2}$, $\forall
i=\overline{1,m}.$ \item[(ii)] If $m\geq 3,$ then
$(a_{ij})_{[k]}=0,$ $\forall k\geq 3,$ $\forall i=\overline{1,m}.$
\end{itemize}
\end{lem}

\begin{pf}  (i).  Let $j\in \alpha_{k_0}$ and
$x=(1-\lambda)e_j+\lambda e_i.$ Here, as before, $e_i, e_j$ are
the vertices of the simplex and $\lambda$ is sufficiently small
positive number. It is easy to see that $x_{[1]}=1-\lambda$ and
$(Vx)_{[1]}=(Vx)_{k_0}.$ Since $Vx\succ x,$ then $x_{[1]}\leq
(Vx)_{[1]},$ so $1-\lambda\leq (Vx)_{k_0}$ or
$$1-\lambda\leq (1-\lambda)^2+2p_{ij,k_0}\lambda(1-\lambda).$$
The last inequality implies that $p_{ij,k_{0}}\geq \frac{1}{2}.$

(ii). Denote $p_{ij,k^\ast}=\max\limits_{t\neq k_{0}}p_{ij,t}.$
One can see that $(Vx)_{k^{\ast}}=(a_{ij})_{[2]}.$ Now from
$$x_{[1]}+x_{[2]}\leq (Vx)_{[1]}+(Vx)_{[2]}$$
 we obtain
 $$1\leq
(1-\lambda)^2+2(p_{ij,k_o}+p_{ij,k^{\ast}})\lambda(1-\lambda).$$
From this inequality we get
$p_{ij,k_o}+p_{ij,k^{\ast}}\geq \frac{2\lambda-{\lambda}^2}{2\lambda(1-\lambda)}=\frac {2-\lambda}{2(1-\lambda)}
=\frac12(1+\frac1{1-\lambda})\ge 1.$ This yields
$p_{ij,k_o}+p_{ij,k^{\ast}}=1$ and $(a_{ij})_{[k]}=0$ $\forall
k\geq 3,$ $\forall i=\overline{1,m}.$
\end{pf}

Now, let us recall some terminology. Let $x^0\in S^{m-1}.$ Then the set $\{x^0,Vx^0,V^2x^0,\cdots\}$ is called the \textit{trajectory } starting at the  point $x^0.$ The point $x^0$  satisfying $Vx^0=x^0$ is called a \textit{fixed}.
The set of all fixed points of the q.s.o. $V$ is denoted by $Fix(V)$ A q.s.o. $V:S^{m-1}\rightarrow S^{m-1}$ is called a {\it regular} if the trajectory of any $x\in S^{m-1}$ converges to
a unique fixed point. Let $V$ be a q.s.o. Then the set $\omega(x^0)=\bigcap \limits_{k
\geq 0}\overline{\bigcup \limits_{n\geq k} \bigl \{ V^nx^0 \bigr
\}}$ is called {\it $\omega$-limit set} of trajectory of initial
point $x^0\in S^{m-1}$. From the compactness of the simplex one
can deduce that $\omega(x^0)\neq \emptyset$ for all $x^0\in
S^{m-1}$.

We provide statements that were proved in \cite{Sh1}.

Let us consider the case $\alpha_1=I$ and $\alpha_k=\emptyset$
for $k\neq 1.$ Then the operator has the following form
$$
\left.
\begin{array}{ll}
(Vx)_1=\sum\limits_{i=1}^mx_i^2+2\sum\limits_{i<j}p_{ij,1}x_ix_j,\\
(Vx)_k=2\sum\limits_{i<j}p_{ij,k}x_ix_j, \ 2\leq k\leq m.\\
\end{array}
\right\} \ \  \eqno (5)
$$

\begin{thm}\label{T1}  A q.s.o. given by (5) is regular. Its
unique fixed point is $e_1.$
\end{thm}

Now we turn to another case, namely let
$\alpha_1=I\backslash\{l\}$ $\alpha_2=\{l\}$  and $\alpha_k=\emptyset, \
\ \forall k\geq 3.$ Then operator (5) has the following form
$$
\left.
\begin{array}{lll}
(Vx)_1=\sum\limits_{i=1, \ i\neq
l}^mx_i^2+2\sum\limits_{i<j}p_{ij,1}x_ix_j,\\[5mm]
(Vx)_2=x_l^2+2\sum\limits_{i<j}p_{ij,2}x_ix_j,\\[5mm]
(Vx)_k=2\sum\limits_{i<j}p_{ij,k}x_ix_j, \ \ 3\leq k\leq m.
\end{array}
\right\} \eqno(6)
$$

\begin{thm}\label{T2} If $l\neq 2$ then the operator (6) is
regular and has a unique fixed point $e_1$. If $l=2$ then the
operator (6) has infinitely many fixed points. Moreover, $\omega$-limit set of trajectory of any
initial point $x^0$ belongs to $co\{e_1,e_2\},$ here $coA$ denotes
the convex hull of a set $A$.
\end{thm}

As we have seen, in the above two theorems the limit behavior of the trajectory of the operator was studied in some particular cases. That is why, it is still important to study the limit behavior of the trajectory of the operator in general case that will be done in the next section.

\section{Main results}

In this section the main results are given. We show that any dissipative q.s.o. has either
a unique or infinitely many fixed points. We also study the limit behavior of the trajectory
of the dissipative q.s.o.

\begin{thm} Any dissipative q.s.o. has either unique or infinitely many fixed points.
\end{thm}

\begin{pf} Since the operator is continuous, and maps convex and compact set into itself,
then according to Bohl-Brouwer theorem there exists a fixed point of the operator.

Further, let us write down a dissipative operator in the following form

$$(Vx)_k=\sum\limits_{i\in
\alpha_k}x_i^2+2\sum\limits_{i<j}p_{ij,k}x_ix_j \ \
k=\overline{1,m},$$ where $\alpha_k\subset
I=\{1,2\cdots,m\},$ $\alpha_i\cap \alpha_j=\emptyset,$ $i\neq j,$
$\bigcup\limits_{k=1}^m\alpha_k=I.$

We prove the statement of the theorem by considering several cases.

\textbf{Case 1.} There is no $k$ such that $k\in \a_k.$ That is, $k\notin \a_k$ for all $k=\overline{1,m}$

\textbf{Case 2.} There exists a unique $k_0$ such that $k_0\in \a_{k_0}.$

\textbf{Case 3.} There are numbers $k_1,k_2,\cdots, k_n (2\le\ n\le m-1)$ such that $k_i\in \a_{k_i} \forall i=\overline{1,n}.$

Now we prove the statement for the Case 1. Let us assume that $k\notin \alpha_k, $ for all $k=\overline{1,m}.$ Then there
exist $k_1,k_2,\cdots,k_n (n<m),$ such that $k_i\in \alpha_{\pi(k_i)},$ where $\pi$ is a
permutation of the set ${k_1,k_2,\cdots, k_n}$ and $\pi(i)\neq i, \ \ \forall i$ due to $k\notin \a_k.$  WLOG we may assume that $k_i=i, \ \forall i=\overline{1,n}$.

Since any permutation is a composition of cycles then $\pi=\pi_1\circ\pi_2\circ\cdots\circ\pi_p,$ where
$\pi_i$ are cycles with the length not less than 2.

Let us consider the case when $p=1.$ WLOG we may assume that $\pi$ is a backward shift, that is, $\pi(i)=i-1$ for $i=\overline{2,n}$ and $\pi(1)=n.$ Then the operator has the following form

$$
\left.
\begin{array}{lll}
(Vx)_k=x_{k+1}^2+\sum\limits_{i\in
\alpha_k\setminus\{k+1\}}x_i^2+\sum\limits_{i<j}2p_{ij,k}x_ix_j, \ \ 1\le k\le n-1\\[5mm]
(Vx)_n= x_1^2+\sum\limits_{i\in
\alpha_n\setminus\{1\}}x_i^2+\sum\limits_{i<j}2p_{ij,n}x_ix_j, \ \\[5mm]
(Vx)_k=\sum\limits_{i\in
\alpha_k}x_i^2+2\sum\limits_{i<j}p_{ij,k}x_ix_j, \ \ n+1\leq k\leq m.
\end{array}
\right\} \eqno(7)
$$

Let $1\le k\le n$ (if $k=n$ then we set $n+1\equiv 1$). Then

$$(Vx)_k=x_{k+1}^2+\sum\limits_{i\in
\alpha_k\setminus\{k+1\}}x_i^2+\sum\limits_{i<j}2p_{ij,k}x_ix_j=$$
$$=x_{k+1}+\sum\limits_{i\in
\alpha_k\setminus\{k+1\}}x_i^2+\sum\limits_{i<j}2p_{ij,k}x_ix_j+x_{k+1}^2-x_{k+1}=$$
$$=x_{k+1}+\sum\limits_{i\in
\alpha_k\setminus\{k+1\}}x_i^2+\sum\limits_{i<j}2p_{ij,k}x_ix_j-\sum\limits_{i=\overline{1,m}, i\neq k+1}x_ix_{k+1}.$$

Thus, the operator (7) can be rewritten as
$$
\left.
\begin{array}{lll}
(Vx)_k=x_{k+1}+\sum\limits_{i\in
\alpha_k\setminus\{k+1\}}x_i^2+\sum\limits_{i<j}2p_{ij,k}x_ix_j-\sum\limits_{i=\overline{1,m}, i\neq k+1}x_ix_{k+1}, \ \ 1\le k\le n-1\\[7mm]
(Vx)_n= x_1+\sum\limits_{i\in
\alpha_n\setminus\{1\}}x_i^2+\sum\limits_{i<j}2p_{ij,n}x_ix_j-\sum\limits_{i=2}^mx_ix_1, \ \\[7mm]
(Vx)_k=\sum\limits_{i\in
\alpha_k}x_i^2+2\sum\limits_{i<j}p_{ij,k}x_ix_j, \ \ n+1\leq k\leq m.
\end{array}
\right\} \eqno(8)
$$

Let us put

$$L_k=\sum\limits_{i<j}2p_{ij,k}x_ix_j-\sum\limits_{i=\overline{1,m}, i\neq k+1}x_ix_{k+1}.\ \eqno (9)$$

We show that $L_k\ge 0,  \ \forall k=\overline{1,n}.$ Indeed, since $k+1\in \alpha_k, \ \forall k=\overline{1,n}$ (with $n+1\equiv 1$) then according to Lemma \ref{l2}, we have $2p_{ik+1,k}\ge 1  \ \forall i=\overline{1,n}, \  i\neq k+1.$Therefore,

$$L_k= 2p_{1k+1,k}x_1x_{k+1}+2p_{2k+1,k}x_2x_{k+1}+\cdots+2p_{kk+1,k}x_kx_{k+1}+2p_{k+2k+1,k}x_{k+2}x_k+\cdots$$
$$+2p_{mk+1,k}x_mx_{k+1}+\sum\limits_{i<j, \ j\neq k+1}2p_{ij,k}x_ix_j-(x_1x_{k+1}+\cdots+x_kx_{k+1}+x_{k+2}x_{k+1}+\cdots+x_mx_{k+1})=$$
$$=(2p_{1k+1,k}-1)x_1x_{k+1}+(2p_{2k+1,k}-1)x_2x_{k+1}+\cdots+(2p_{kk+1,k}-1)x_kx_{k+1}+$$$$+(2p_{k+2k+1,k}-1)x_{k+2}x_k+\cdots+(2p_{mk+1,k}-1)x_mx_{k+1}+\sum\limits_{i<j, \ j\neq k+1}2p_{ij,k}x_ix_j\ge 0.$$

Further, since $\pi$ is a permutation with maximal length, then

$$\sum\limits_{i\in \a_k}x_i^2=0 \ \ \eqno (10)$$
for $n+1\le k\le m.$ Indeed, the set $\a_k$ for $n+1\le k\le m$ does not contain numbers $1,2,\cdots,n$ because these numbers belong only to one of the sets $\a_1,\a_2,\cdots,\a_n.$ Moreover, $k\notin \a_k$ according to the assumption. If $j_0\in \a_k$ for $j_0\ge n+1 \ j\neq k+1$ then there exist numbers
$j_i,j_2,\cdots,j_r (j_r\ge n+1, \ r\ge 1),$ such that $j_0\in\a_k,\ j_1\in\a_{j_0},\cdots j_r\in\a_{j_{r-1}},\ k\in\a_{j_r}.$ Let $\pi'$ be a permutation such that $\pi'(k)=j_0,\ \pi'(j_0)=j_1,\cdots, \pi'(j_{r-1})=j_r,\ \pi'(j_r)=k,$ then $\sigma=\pi\circ\pi'$ is the permutation on $\{1,2,\cdots,m\}$  with greater length than $\pi,$ which is a contradiction. That is why, $\a_k=\emptyset$ for $n+1\le k\le m.$ Hence, numbers $k+1,k+2,\cdots, m$ belong to one of the sets $\a_1,\a_2,\cdots,\a_n.$

Since $L_k\ge 0$ for all $k=\overline{1,n},$ then by letting $Vx=x$ and by summing up (8) from $k=1$ to $n$ we get

$$\sum\limits_{k=1}^n\sum\limits_{i\in\a_k\setminus\{k+1\}}x_i^2=0.$$

Due to $(4')$, it is clear that $\sum\limits_{k=1}^m\sum\limits_{i\in \a_k}x_i^2=\sum\limits_{k=1}^mx_i^2.$ Moreover, from (10) we get $\sum\limits_{k=1}^n\sum\limits_{i\in \a_k}x_i^2=\sum\limits_{k=1}^mx_i^2.$
 Therefore,
$$0=\sum\limits_{k=1}^n\sum\limits_{i\in\a_k\setminus\{k+1\}}x_i^2=\sum\limits_{k=1}^n\sum\limits_{i\in\a_k}x_i^2-
\sum\limits_{k=1}^nx_i^2=\sum\limits_{k=1}^mx_i^2-\sum\limits_{k=1}^nx_i^2=\sum\limits_{i=n+1}^mx_i^2,$$
which implies $x_{n+1}=x_{n+2}=\cdots=x_{m}=0.$

In addition, we obtain $x_1=x_2=\cdots=x_n$ which means that
$(\frac1n,\frac1n,\cdots,\frac1n,0,0\cdots,0)$ is a unique fixed point.

Now, let $p\ge 2.$ There is no loss in generality in assuming $p=2,$ that is $\pi=\pi_1\circ\pi_2.$
WLOG we may put $\pi_1(i)=i-1$ for $i=\overline{2,n_1}$ with $\pi_1(1)=n_1$ and $\pi_2(i)=i-1$ for $i=\overline{n_1+2,n_2}$ with $\pi_2(n_1+1)=n_2.$ Then the operator has the following form.

$$
\left.
\begin{array}{lllll}
  (Vx)_k=x_{k+1}^2+\sum\limits_{i\in
\alpha_k\setminus\{k+1\}}x_i^2+\sum\limits_{i<j}2p_{ij,k}x_ix_j, \ \ 1\le k\le n_1-1\\[5mm]
  (Vx)_{n_1}= x_1^2+\sum\limits_{i\in
\alpha_{n_1}\setminus\{1\}}x_i^2+\sum\limits_{i<j}2p_{ij,n_1}x_ix_j, \ \\[5mm]
   (Vx)_k=x_{k+1}^2+\sum\limits_{i\in
\alpha_k\setminus\{k+1\}}x_i^2+\sum\limits_{i<j}2p_{ij,k}x_ix_j, \ \ n_1+1\le k\le n_2-1\\[5mm]\\
  (Vx)_{n_2}= x_{n_1+1}^2+\sum\limits_{i\in
\alpha_{n_2}\setminus\{n_1+1\}}x_i^2+\sum\limits_{i<j}2p_{ij,n_2}x_ix_j, \ \\[5mm] \\
  (Vx)_k=\sum\limits_{i\in
\alpha_k}x_i^2+2\sum\limits_{i<j}p_{ij,k}x_ix_j, \ \ n_2+1\leq k\leq m
\end{array}
\right\} \eqno (11)
$$

which can be rewritten as

$$
\left.
\begin{array}{lllll}
  (Vx)_k=x_{k+1}+\sum\limits_{i\in
\alpha_k\setminus\{k+1\}}x_i^2+\sum\limits_{i<j}2p_{ij,k}x_ix_j-\sum\limits_{i=\overline{1,m}, i\neq k+1}x_ix_{k+1}, \ \ 1\le k\le n_1-1\\[5mm]
  (Vx)_{n_1}= x_1+\sum\limits_{i\in
\alpha_{n_1}\setminus\{1\}}x_i^2+\sum\limits_{i<j}2p_{ij,n_1}x_ix_j-\sum\limits_{i=2}^mx_ix_1, \ \\[5mm]
   (Vx)_k=x_{k+1}+\sum\limits_{i\in
\alpha_k\setminus\{k+1\}}x_i^2+\sum\limits_{i<j}2p_{ij,k}x_ix_j-\sum\limits_{i=\overline{1,m}, i\neq k+1}x_ix_{k+1}, \ \ n_1+1\le k\le n_2-1\\[5mm]\\
  (Vx)_{n_2}= x_{n_1+1}+\sum\limits_{i\in
\alpha_{n_2}\setminus\{n_1+1\}}x_i^2+\sum\limits_{i<j}2p_{ij,n_2}x_ix_j-\sum\limits_{i=\overline{1,m}, i\neq n_1+1}x_ix_{n_1+1}, \ \\[5mm] \\
  (Vx)_k=\sum\limits_{i\in
\alpha_k}x_i^2+2\sum\limits_{i<j}p_{ij,k}x_ix_j, \ \ n_2+1\leq k\leq m.
\end{array}
\right\} \eqno (12)
$$

If we put $Vx=x$ and sum up (12) from $k=1$ to $k=n_2,$ we obtain

$$\sum\limits_{k=1}^{n_2}(Vx)_k=\sum\limits_{k=1}^{n_2}x_k+S_1+S_2+\sum\limits_{k=1}^{n_2}L_k,$$
where $S_1=\sum\limits_{k=1}^{n_1}\sum\limits_{i\in \a_k\setminus\{k+1\}}x_i^2 \ (n_1+1\equiv 1)$
and $S_2=\sum\limits_{k=n_1+1}^{n_2}\sum\limits_{i\in \a_k\setminus\{k+1\}}x_i^2 \ (n_2+1\equiv n_1+1).$
Since $L_k\ge 0\ \ \forall k=\overline{1, n_2}$ and $\sum\limits_{k=1}^{n_2}(Vx)_k=\sum\limits_{k=1}^{n_2}x_k,$ then we have
$$\sum\limits_{k=1}^{n_2}x_k=\sum\limits_{k=1}^{n_2}(Vx)_k=\sum\limits_{k=1}^{n_2}x_k+S_1+S_2+\sum\limits_{k=1}^{n_2}L_k\ge
\sum\limits_{k=1}^{n_2}x_k+ S_1+S_2 $$
which implies $S_1+S_2=0.$ As we showed the relation (10), we can also show that $\sum\limits_{i\in \a_k}x_i^2$
for $n_2+1\le k\le m.$ Therefore, $0=S_1+S_2=\sum\limits_{k=n_2+1}^mx_i^2,$ and hence $x_i=0\ \ \forall i\ge n_2+1.$

Now, since $\sum\limits_{i\in \a_k\setminus\{k+1\}}x_i^2$ and $L_k$ are positive, we let $Vx=x,$ and from (12) we get
$x_k\ge x_{k+1}$ for $k=\overline{1,n_1} \ (x_{n_1+1}\equiv x_1),$ which implies $x_1=x_2=\cdots=x_{n_1}=\a.$
Consequently, we get $x_k\ge x_{k+1}$ for $k=\overline{n_1+1,n_2} \ (x_{n_2+1}\equiv x_{n_1+1}),$ which implies $x_{n_1+1}=x_{n_1+2}=\cdots=x_{n_2}=\beta.$ So, the set of fixed points is

$$\{(\a,\a,\cdots,\a,\beta,\beta,\cdots,\beta,0,0,\cdots,0)\in S^{m-1}|\ \ n_1\a+n_2\beta=1\}$$
Thus, there are infinitely many fixed points.

\textbf{Case 2.}
Assume that there exist a unique $k_1\in \alpha_{k_1},$ and WLOG we assume that $k_1=1$
There are two possibilities

(a). There is a permutation $\pi$ of the numbers $k_2,k_3,\cdots,k_n\ (2\le k_i\le m,\ 1\le n\le m-1)$ such that $k_i\in \a_{\pi(k_i)}.$

(b). There is no any permutation $\pi$ of the set such that $k_i\in \a_{\pi(k_i)}.$

In case (a), WLOG we may assume that $k_i=i$ and $\pi(i)=i-1$ for $i=\overline{3,n}$ and $\pi(2)=n$ then the operator has the following form

$$
\left.
\begin{array}{llll}
  (Vx)_1=x_1^2+\sum\limits_{i\in
\alpha_1\setminus\{1\}}x_i^2+\sum\limits_{i<j}2p_{ij,1}x_ix_j, \ \ \\[5mm]
   (Vx)_k=x_{k+1}^2+\sum\limits_{i\in
\alpha_k\setminus\{k+2\}}x_i^2+\sum\limits_{i<j}2p_{ij,k}x_ix_j, \ \ 2\le k\le n-1\\[5mm]\\
   (Vx)_n= x_2^2+\sum\limits_{i\in
\alpha_n\setminus\{2\}}x_i^2+\sum\limits_{i<j}2p_{ij,n}x_ix_j, \ \\[5mm]
    (Vx)_k=\sum\limits_{i\in
\alpha_k}x_i^2+2\sum\limits_{i<j}p_{ij,k}x_ix_j, \ \ n+1\leq k\leq m,
\end{array}
\right\} \eqno (13)
$$

which can be rewritten as

$$
\left.
\begin{array}{llll}
  (Vx)_1=x_1+\sum\limits_{i\in
\alpha_1\setminus\{1\}}x_i^2+\sum\limits_{i<j}2p_{ij,1}x_ix_j-\sum\limits_{i=2}^mx_ix_1, \ \ \\[5mm]
   (Vx)_k=x_{k+1}+\sum\limits_{i\in
\alpha_k\setminus\{k+2\}}x_i^2+\sum\limits_{i<j}2p_{ij,k}x_ix_j-\sum\limits_{i=\overline{1,m}, i\neq k+2}x_ix_{k+2}, \ \ 2\le k\le n-1\\[5mm]\\
   (Vx)_n= x_2+\sum\limits_{i\in
\alpha_n\setminus\{2\}}x_i^2+\sum\limits_{i<j}2p_{ij,n}x_ix_j-\sum\limits_{i=\overline{1,m}, i\neq 2}x_ix_2, \ \\[5mm]
    (Vx)_k=\sum\limits_{i\in
\alpha_k}x_i^2+2\sum\limits_{i<j}p_{ij,k}x_ix_j, \ \ n+1\leq k\leq m.
\end{array}
\right\} \eqno (14)
$$

We use the similar arguments as in Case 1 and deduce that the set of fixed points of (14) is
$$\{(1-(n-1)\a,\underbrace{\a,\a,\cdots,\a}_{n-1\  times},0,0\cdots 0)| \ \ 0<\a<\frac1n \}$$

Let us now assume that (b) holds. The operator has the following form
$$
\left.
\begin{array}{ll}
  (Vx)_1=x_1^2+\sum\limits_{i\in
\a_1\setminus\{1\}}x_i^2+\sum\limits_{i<j}2p_{ij,1}x_ix_j, \ \ \\[5mm]
  (Vx)_k=\sum\limits_{i\in
\a_k}x_i^2+2\sum\limits_{i<j}p_{ij,k}x_ix_j, \ \ 2\leq k\leq m,
\end{array}
\right\} \eqno (15)
$$
 which can be rewritten as
$$
\left.
\begin{array}{ll}
  (Vx)_1=x_1+\sum\limits_{i\in
\a_1\setminus\{1\}}x_i^2+\sum\limits_{i<j}2p_{ij,1}x_ix_j-\sum\limits_{i=2}^mx_ix_1, \ \ \\[5mm]
  (Vx)_k=\sum\limits_{i\in
\a_k}x_i^2+2\sum\limits_{i<j}p_{ij,k}x_ix_j, \ \ 2\leq k\leq m.
\end{array}
\right\} \eqno (16)
$$

If (b) holds, then it is easy to show (as was done in (10)) that
$$\sum\limits_{i\in \a_k}x_i^2=0,\ \ \ \forall k=\overline{2,m}.$$

Hence, $\sum\limits_{i\in\a_1\setminus\{1\}}x_i^2=\sum\limits_{i=2}^mx_i^2.$ So, if we set $Vx=x $ in (16), then it follows that $x_2=x_3=\cdots=x_m=0,$ which means that $e_1$ is a unique fixed point for (16).

\textbf{Case 3.} Let us suppose that there are numbers $k_1,k_2,\cdots, k_n\ (2\le\ n\le m-1)$ such that $k_i\in \a_{k_i}, \forall i=\overline{1,n}.$ It is sufficient to show for $n=2$ and  $k_1=1, \ k_2=2.$

Then the operator has the following form

$$
\left.
\begin{array}{lll}
  (Vx)_1=x_1^2+\sum\limits_{i\in
\a_1\setminus\{1\}}x_i^2+\sum\limits_{i<j}2p_{ij,1}x_ix_j, \ \ \\[5mm]
  (Vx)_2=x_2^2+\sum\limits_{i\in
\a_2\setminus\{2\}}x_i^2+\sum\limits_{i<j}2p_{ij,2}x_ix_j, \ \ \\[5mm]
  (Vx)_k=\sum\limits_{i\in
\a_k}x_i^2+2\sum\limits_{i<j}p_{ij,k}x_ix_j, \ \ 3\leq k\leq m,
\end{array}
\right\} \eqno (17)
$$
 which can be rewritten as
$$
\left.
\begin{array}{lll}
  (Vx)_1=x_1+\sum\limits_{i\in
\a_1\setminus\{1\}}x_i^2+\sum\limits_{i<j}2p_{ij,1}x_ix_j-\sum\limits_{i=2}^mx_ix_1, \ \ \\[5mm]
   (Vx)_2=x_2+\sum\limits_{i\in
\a_2\setminus\{2\}}x_i^2+\sum\limits_{i<j}2p_{ij,2}x_ix_j-\sum\limits_{i=\overline{1,m},\ i\neq 2}x_ix_1, \ \ \\[5mm]
  (Vx)_k=\sum\limits_{i\in
\a_k}x_i^2+2\sum\limits_{i<j}p_{ij,k}x_ix_j, \ \ 3\leq k\leq m.
\end{array}
\right\} \eqno (18)
$$

Now, if there exist numbers $l_3,l_4,\cdots,l_n\ (3\le l_i\le m,\ 1\le n\le m-2)$ such that $l_i\in \a_{\pi(l_i)},$ then by putting $l_i=i$ and $\pi(i)=i-1$ with $\pi(3)=n$ we find that the set of all fixed points is
$$\{(\a,\beta, \underbrace{\gamma,\gamma,\cdots,\gamma}_{n-2\  times},0,0\cdots 0)| \ \ \a+\beta+(n-2)\gamma=1 \}.$$

If there is no such numbers then by using arguments as (10) was shown, we can also show that
$$\sum\limits_{i\in \a_k}x_i^2=0,\ \ \ \forall k=\overline{3,m},$$
which means that $\sum\limits_{i\in\a_1\setminus\{1\}}x_i^2+\sum\limits_{i\in
\a_2\setminus\{2\}}x_i^2=\sum\limits_{i=3}^mx_i^2.$ So, if we put $Vx=x $ in (18), then it follows that $x_3=x_4=\cdots=x_m=0,$ which means that the vectors $\lambda e_1+(1-\lambda)e_2$ are fixed points of (18).
\end{pf}

\begin{rem} It can be seen from the proof of the theorem that the set of fixed points of operator has shown exactly
up to permutation of coordinates. In those cases when operator has a unique fixed point, this fixed point is
either the vertex of a simplex or the center of a face of the simplex. In those cases when operator has
infinitely many fixed points, then the set of fixed points is some face of the simplex.
\end{rem}

Now let us consider some examples.

\textbf{Example 1.} Consider the following operator on 2D simplex.

\begin{eqnarray*}
&&(Vx)_1=x_2x_3,\\
&&(Vx)_2=x_1^2+x_1x_2+x_1x_3,\\
&&(Vx)_3=x_2^2+x_3^2+x_1x_2+x_1x_3+x_2x_3.
\end{eqnarray*}
By putting $Vx=x$ and solving the system of equations we find that $(0,0,1)$ is a unique fixed point.

\textbf{Example 2.} Consider the following operator on 3D simplex.

\begin{eqnarray*}
&&(Vx)_1=x_1^2+x_1x_2+x_1x_3+ax_1x_4+x_2x_4,\\
&&(Vx)_2=x_3^2+x_1x_3+x_2x_3+bx_3x_4,\\
&&(Vx)_3=x_2^2+x_1x_2+x_2x_3,\\
&&(Vx)_4= (2-a)x_1x_4+(2-b)x_3x_4,
\end{eqnarray*}

where $1\le a,b\le 2.$ The set of fixed points of the above operator is $(1-2\lambda,\lambda,\lambda,0),$ where $0\le\lambda\le\frac12.$

\textbf{Example 3.} Consider the following operator on 2D simplex.

\begin{eqnarray*}
&&(Vx)_1=x_2^2+x_3^2+x_1x_2+x_1x_3+x_2x_3 ,\\
&&(Vx)_2=x_2x_3,\\
&&(Vx)_3=x_1^2+x_1x_2+x_1x_3.
\end{eqnarray*}

Above operator has a unique point $(\frac12,0,\frac12).$

In general, for any given vertex of the simplex or any given center of a face of the simplex, it is possible to provide an example of dissipative q.s.o. with a unique fixed point at the given point.

\begin{thm} Let $V$ be dissipative q.s.o. Then one of the following statements always holds.
\begin{itemize}
  \item The operator is regular. Its unique point is either a vertex of the simplex or the center of its face.
  \item The operator has infinitely many fixed points. $\omega-$ limit set of any initial point is contained in the set of fixed points i.e. $\omega(x)\subset Fix(V)$.
\end{itemize}
\end{thm}

\begin{pf}
The way of the arguments in the proof is similar to the proof of the previous theorem. Here we divide the proof into 3 cases as well.

\textbf{Case 1.} In this case the operator has one of the forms (8) and (11).
Suppose that the operator is of form (8).
Define
$$\varphi(x)=x_1+x_2+\cdots+x_n.$$

Then
$$\varphi(Vx)=\sum\limits_{i=1}^n(Vx)_i=\sum\limits_{i=1}^nx_i+\sum\limits_{k=1}^n\sum\limits_{i\in\a_k\setminus\{k+1\}}x_i^2
+\sum\limits_{i=1}^nL_i. \ \eqno (19)$$

Here $L_i$ satisfy (9). Since $L_i\ge 0,$ then from (19) it follows that
$\varphi(Vx)\ge\varphi(x).$ Hence, the sequence $\{\varphi(V^k(x)): k=1,2,\cdots\}$ is non-decreasing and bounded. That is why the limit $\lim\limits_{k\to\infty}\varphi(V^kx)$ exists.
Let us put $C=\lim\limits_{k\to\infty}\varphi(V^kx).$

Taking (10) into account one can rewrite (19) as

$$\varphi(Vx)=\sum\limits_{i=1}^n(Vx)_i=\sum\limits_{i=1}^nx_i+\sum\limits_{i=n+1}^mx_i^2
+\sum\limits_{i=1}^nL_i. \ \eqno (20)$$

From (20) we get

$$\varphi(V^{k+1}x)=\sum\limits_{i=1}^n(V^{k+1}x)_i=\sum\limits_{i=1}^n(V^kx)_i+\sum\limits_{i=n+1}^m(V^kx)_i^2
+\sum\limits_{i=1}^nL^k_i=$$
$$=\varphi(V^{k}x)+ \sum\limits_{i=n+1}^m(V^kx)_i^2+\sum\limits_{i=1}^nL^k_i,\ \eqno (21)$$
where

$$L^k_s=\sum\limits_{i<j}2p_{ij,s}(V^kx)_i(V^kx)_j-\sum\limits_{i=\overline{1,m}, i\neq s+1}(V^kx)_i(V^kx)_{s+1}\ \ s=\overline{1,n},\ \ k=1,2,\cdots.$$

Since $\lim\limits_{k\to\infty}\varphi(V^{k+1}x)=\lim\limits_{k\to\infty}\varphi(V^{k}x)=C$ and $L^k_s\ge 0,$ then (21) implies $\sum\limits_{i=n+1}^m(V^kx)_i^2\to 0$ as $k\to\infty,$ Therefore, $(V^kx)_i\rightarrow 0$ for all $i=\overline{n+1,m},$ as $k\to\infty.$  Taking into account $\sum\limits_{i=1}^m(V^kx)_i=1$  we get $C=1.$ Furthermore, if $k\to\infty$ then $\sum\limits_{i\in\a_s\setminus\{s+1\}}(V^kx)_i^2\to 0$ for $s=\overline{1,n}$ and $L^k_i\to 0.$ Therefore, from

$$(V^{k+1}x)_i=(V^kx)_{i+1}+\sum\limits_{i\in\a_i\setminus\{i+1\}}(V^kx)_i^2+L^n_i, \ \ 1\le i\le n,$$
we get
$$\lim\limits_{k\rightarrow\infty}(V^kx)_{i+1}=\lim\limits_{k\rightarrow\infty}(V^kx)_i,$$
which implies $\lim\limits_{k\rightarrow\infty}(V^kx)_i=\frac1n, \ \forall i=\overline{1,n}.$ This means that the operator is regular with the unique fixed point $(\frac1n,\frac1n,\cdots,\frac1n,0,0,\cdots,0).$

Now, let us assume that the operator is of form (11). Using the above method, we define

$$\varphi(x)=x_1+x_2+\cdots+x_{n_2}$$
and show that $\lim\limits_{k\to\infty}\varphi(V^kx)=1.$ Moreover, it can also be shown that
$$\lim\limits_{k\rightarrow\infty}(V^kx)_{i+1}=\lim\limits_{k\rightarrow\infty}(V^kx)_i$$
for $i=\overline{1,n_1}$ and
$$\lim\limits_{k\rightarrow\infty}(V^kx)_{i+1}=\lim\limits_{k\rightarrow\infty}(V^kx)_i$$
for $i=\overline{n_1+1,n_2}.$ The last equalities imply that the trajectory of any initial point approaches the set of fixed points of the operator.

Now, we turn to Case 2. Assume that (a) holds. Then we define
$$\varphi(x)=x_1+x_2+\cdots+x_{n}$$
and show that $\lim\limits_{k\to\infty}\varphi(V^kx)=1.$ Its proof was done in Case 1, so we do not repeat it.
Moreover,  we have

$$\lim\limits_{k\rightarrow\infty}(V^kx)_{i+1}=\lim\limits_{k\rightarrow\infty}(V^kx)_i$$
for $i=\overline{2,n-1}.$ That is why, for any initial point $x^0\in S^{m-1}$ we have

$$\omega(x^0)\subset \{(\a,\beta,\cdots,\beta,0,\cdots,0)|\ \ \a+(n-1)\beta=1\}=Fix(V).$$

Now, if (b) holds, then we put $\varphi(x)=x_1$ and show that $\lim\limits_{k\to\infty}\varphi(V^kx)=1,$ which means that V is regular at $e_1.$

Let us consider the last Case 3 in which $V$ is of form (18). If there do not exist numbers $l_3,l_4,\cdots,l_n\ \ (3\le l_i\le m,\ 1\le n\le m-2)$ such that $l_i\in \a_{\pi(l_i)}$ for some permutation $\pi,$  then we define
$\varphi(x)=x_1+x_2$ and show that $\lim\limits_{k\to\infty}\varphi(V^kx)=1,$ which means that $\omega(x)\subset Fix(V).$

Otherwise we put $l_i=i$ and $\pi(i)=i-1$ with $\pi(3)=n.$ Then we define
$$\varphi(x)=x_1+x_2+\cdots+x_{n}$$
and show that $\lim\limits_{k\to\infty}\varphi(V^kx)=1.$ Furthermore,

$$\lim\limits_{k\rightarrow\infty}(V^kx)_{i+1}=\lim\limits_{k\rightarrow\infty}(V^kx)_i$$ holds for
$i=\overline{3,n-1}.$ Hence,

$$\omega(x^0)\subset \{(\a,\beta,\gamma,\cdots,\gamma,0,\cdots,0)|\ \ \a+\beta+(n-2)\gamma=1\}=Fix(V).$$

\end{pf}

\begin{cor}Let $V$ be a dissipative q.s.o. on 2D simplex.

If $V$ is regular, then its unique fixed point is one of the following points:
$$e_1, e_2, e_3, \frac12(e_1+e_2), \frac12(e_1+e_3), \frac12(e_2+e_3).$$

If $V$ has infinitely many fixed points, then either $\omega(x)\subset \{\lambda e_1+(1-\lambda)e_2|\ 0\le\lambda\le 1\}$
or $\omega(x)\subset \{\lambda e_1+(1-\lambda)e_3|\ 0\le\lambda\le 1\}$ or $\omega(x)\subset \{\lambda e_2+(1-\lambda)e_3|\ 0\le\lambda\le 1\}.$
\end{cor}

\section{Conclusion}

The main achievement of this paper is the classification of the limit behavior of the trajectory of dissipative q.s.o. on finite dimensional simplex. The main results are Theorem 3.1, Theorem 3.3 and Corollary 3.4. Most of the papers on quadratic stochastic operators are considered in small dimensions of the simplex. But, in the present paper it was considered in any finite dimensional simplex which made the investigation more difficult. The methods, which were used for proving results can be used in other disciplines of mathematics, namely nonlinear analysis, dynamical systems and ergodic theory. Here some directions in which research can be carried on.

1. One can consider polynomial or any continuous function $V,$ defined on the simplex satisfying $Vx\succ x$ and study the limit behavior of its trajectory. Research can be started when $V$ is a cubic operator, that is $V(\lambda x)=\lambda^3Vx, \ \lambda\ge 0.$

2. Dissipative quadratic stochastic operators can be defined on infinite dimensional space. Infinite dimensional simplex can be defined on $\ell_1$ as $S=\{x\in\ell_1: x_{i}\geq 0, \ \ \sum\limits_{i=1}^{\infty}x_{i}=1\}.$ Majorization for infinite cases is developed. But, the problem is still difficult since infinite dimensional simplex is not compact. We also can not state that dissipative q.s.o. has at least one fixed point.

\section{Acknowledgement}
The authors wish to thank the referees for their valuable comments and suggestions, which improved the quality of the paper.

\end{document}